\documentclass[11pt]{amsart}
\usepackage[english]{babel}
\usepackage{graphicx}
\newtheorem{theorem}{Theorem}[section]
\newtheorem{corollary}[theorem]{Corollary}
\newtheorem{lemma}[theorem]{Lemma}

\theoremstyle{definition}

\theoremstyle{remark}
\newtheorem{remark}{Remark}
\numberwithin{equation}{section}
\newcommand{\real}{\mathbb R}
\def\natu{\mathbb N}
\def\dis{\displaystyle}
\def\sobo1{H^{1}(0,L)}
\def\sobol{(H^{1}(0,L))^n}

\def\intl{\int_{0}^{L}}
\def\intomega{\int_{\Omega}}

\def\comega{\overline{\Omega}}

\def\afirm#1#2#3{\skipaline\noindent{\bf #1} \hfill
\begin{tabular}{p{#2}}{\sl #3} \end{tabular}\hfill $ $\skipaline}
\def\skipaline{\vskip12pt plus 1pt}

\begin{document}

\title[Matrix Lyapunov inequalities]
{Matrix Lyapunov inequalities for ordinary and elliptic partial differential equations}%
\author{Antonio Ca\~{n}ada}
\author{Salvador Villegas}
\thanks{The authors have been supported by the Ministry of Education and
Science of Spain (MTM2005-01331)}
\address{Departamento de An\'{a}lisis
Matem\'{a}tico, Universidad de Granada, 18071 Granada, Spain.}
\email{acanada@ugr.es, svillega@ugr.es}

\subjclass[2000]{34B15, 34B05, 35J20, 35J25, 35J65}%
\keywords{Neumann boundary value problems, Matrix Lyapunov
inequalities, ordinary differential equations, elliptic partial differential equations, resonant problems.}%

\begin{abstract}
This paper is devoted to the study of $L_p$ Lyapunov-type
inequalities for linear systems of equations with Neumann boundary
conditions and for any constant $p \geq 1$. We consider ordinary
and elliptic problems. The results obtained in the linear case are
combined with Schauder fixed point theorem to provide new results
about the existence and uniqueness of solutions for resonant
nonlinear problems. The proof uses in a fundamental way the
nontrivial relation between the best Lyapunov constants and the
minimum value of some especial minimization problems.
\end{abstract}
\maketitle
\section{Introduction}

\noindent Let us consider the linear Neumann boundary problem
\begin{equation}\label{n1}
u''(x) + a(x)u(x) = 0, \ x \in (0,L), \ u'(0) = u'(L) = 0
\end{equation}
and let $1 \leq p \leq \infty$ be given. If function $a$ satisfies
\begin{equation}\label{2806072}
a \in L^p(0,L)\setminus \{0\},\  \displaystyle \int_{0}^{L} a(x) \
dx \geq 0, \end{equation} $L_p-$Lyapunov inequality provides
optimal necessary conditions for  boundary value problem
(\ref{n1}) to have nontrivial solutions, given in terms of the
$L^p$ norm, $\Vert \cdot \Vert_p,$ of the function $a^+$, where
$a^{+}(x) = \max \{ a(x),0 \}$(see \cite{ha} and \cite{huyo} for
the case $p=1$ and \cite{camovi}, \cite{zhang} for the case $1 < p
\leq \infty$).

In particular, under the restriction (\ref{2806072}) for $p=1,$
$L_1-$Lyapunov inequality may be used to prove that (\ref{n1}) has
only the trivial solution if function $a$ satisfies
\begin{equation}\label{c1}
\int_{0}^{L} a^{+}(x) \ dx \leq 4/L
\end{equation}
In a similar way, under (\ref{2806072}) for $p= \infty,$
$L_\infty-$Lyapunov inequality may be used to prove that
(\ref{n1}) has only the trivial solution if function $a$ satisfies
\begin{equation}\label{c2}
a^+ \prec \pi^2/L^2,
\end{equation}
where for $c,d \in L^1 (0,L),$ we write $c \prec d$ if $c(x) \leq
d(x)$ for a.e. $x \in [0,L]$ and $c(x) < d(x)$ on a set of
positive measure. Moreover, (\ref{c1}) and (\ref{c2}) are,
respectively, optimal $L_1$ and $L_\infty$ restrictions (see
Remark \ref{2409081} below).

If $p=\infty,$ assumptions (\ref{2806072}) and (\ref{c2}) are a
nonuniform nonresonance condition with respect to the two first
eigenvalues $\lambda_0 = 0 $ and $\lambda_1 = \pi^2/L^2$ of the
eigenvalue problem
\begin{equation}\label{n2}
 u''(x) + \lambda u(x) = 0, \ x \in (0,L), \ u'(0) = u'(L) =
0 \end{equation} (see \cite{mawa}) while if $p=1,$ (\ref{c1}) was
first introduced by Lyapunov under Dirichlet boundary conditions
(see \cite{ha}, chapter XI, for some generalizations and historic
references and \cite{cavi2} for $L_1-$Lyapunov inequality at
higher eigenvalues).

It is clear that (\ref{c1}) and (\ref{c2}) are not related. A
natural link between them arises if $L_p-$Lyapunov inequalities,
for $1 < p < \infty,$ are considered and then one examines what
happens if $p \rightarrow 1^{+}$ and $p \rightarrow \infty$
(\cite{camovi}). One of the main applications of Lyapunov
inequalities is its use in the study of nonlinear resonant
problems.

Different authors have generalized the $L_\infty-$Lyapunov
inequality (\ref{2806072})-(\ref{c2}) to vector differential
equations of the form
\begin{equation}\label{s0}
u''(x) + A(x)u(x) = 0, \ x \in (0,L)
\end{equation}
where $A(\cdot)$ is a real and continuous $n\times n$ symmetric
matrix valued function, together with different boundary
conditions. These $L_\infty$ generalizations have been given not
only at the two first eigenvalues but also at higher eigenvalues
of (\ref{n2}) and they have been used in the study of resonant
nonlinear problems (\cite{ah}, \cite{brli}, \cite{la},
\cite{lasa}, \cite{wa}). Also, some abstract versions for
semilinear equations in Hilbert spaces and applications to
elliptic problems and semilinear wave equations have been given in
\cite{ba}, \cite{foma}, \cite{ma} and \cite{ma2}. In spite of its
interest in the study of different questions such as stability
theory, the calculation of lower bounds on eigenvalue problems,
etc. (\cite{clhi}, \cite{ha}, \cite{zhang}), the use of
$L_\infty-$Lyapunov inequalities in the study of nonlinear
resonant problems only allows a weak interaction between the
nonlinear term and the spectrum of the linear part. For example,
using the $L_\infty-$Lyapunov inequalities showed in \cite{la} for
the periodic boundary value problem (see also \cite{ah} and
\cite{brli}), it may be proved that if there exist real symmetric
matrices $P$ and $Q$ with eigenvalues $p_1 \leq \cdots \leq p_n$
and $q_1 \leq \cdots \leq q_n,$ respectively, such that
\begin{equation}\label{eq2307081}
P \leq G''(u) \leq Q, \ \forall \ u \in \real^n
\end{equation}
and such that
\begin{equation}\label{eq2307082}
\bigcup_{i=1}^{n} \ [p_i,q_i] \cap \{ k^2: k \in \natu \cup \{ 0
\} \} = \emptyset,
\end{equation}
then, for each continuous and $2\pi-$periodic function $h$, the
periodic problem
\begin{equation}\label{eq2307083}
u''(x) + G'(u(x)) = h(x), x \in (0,2\pi), u(0) - u(2\pi) = u'(0) -
u'(2\pi) = 0,
\end{equation}
has a unique solution. Here $G: \real^n \rightarrow \real$ is a
$C^2-$mapping and the relation $C \leq D$ between $n\times n$
matrices means that $D-C$ is positive semi-definite. Now, by using
the variational characterization of the eigenvalues of a real
symmetric matrix, it may be easily deduced that (\ref{eq2307081})
and (\ref{eq2307082}) imply that the eigenvalues $g_1 (u) \leq
\cdots \leq g_n (u)$ of the matrix $G''(u),$ satisfy
\begin{equation}\label{eq2307084}
p_i \leq g_i (u) \leq q_i, \ \forall \ u \in \real^n.
\end{equation}
Consequently each continuous function $g_i (u), \ 1 \leq i \leq
n,$ must fulfil
\begin{equation}\label{eq2307085}
g_i (\real^n) \cap \{ k^2: k \in \natu \cup \{ 0 \} \} =
\emptyset.
\end{equation}
To the best of our knowledge, we do not know any previous work on
$L_p$ Lyapunov inequalities when $1 \leq p < \infty$ for systems
of the type (\ref{s0}) under Neumann boundary conditions. Really,
if the restrictions on the matrix $A(x)$ are of $L_p$ type, with
$1 \leq p < \infty,$ it seems difficult to use the ideas contained
in the mentioned papers to get new results on problems at
resonance.

In the second section of this paper we provide for each $p,$ with
$1 \leq p \leq \infty,$ optimal necessary conditions for boundary
value problem
\begin{equation}\label{s1}
u''(x) + A(x)u(x) = 0, \ x \in (0,L), \ u'(0) = u'(L) = 0,
\end{equation}
to have nontrivial solutions. These conditions are given in terms
of the $L^p$ norm of appropriate functions $b_{ii}(x), \ 1 \leq i
\leq n,$ related to $A(x)$ through the inequality $A(x) \leq B(x),
\ \forall\ x \in [0,L],$ where $B(x)$ is a diagonal matrix with
entries given by $b_{ii}(x), \ 1 \leq i \leq n.$ In particular, we
can use different $L_{p_i}$ criteria for each $1 \leq i \leq n$
and this confers a great generality on our results. Even in the
case $p_i = \infty,\ 1 \leq i \leq n,$ our method of proof is
different from those given in previous works. In fact, we begin
Section 2 with a lemma inspired from \cite{la} and \cite{lasa},
where the authors studied the periodic problem. The proof that we
give for this lemma suggest the way for the case when $1 \leq p <
\infty,$ where we use in a fundamental way some previous results
which have been proved in \cite{camovi} and \cite{camovi2}. They
relate, for ordinary and elliptic problems, the best Lyapunov
constants to the minimum value of some especial minimization
problems. If $1 < p < \infty,$ this minimum value plays the same
role as, respectively, the constants $4/L$ (if $p=1$)  and $\pi^2
/L^2$ (if $p=\infty$) in (\ref{c1}) and (\ref{c2}) (see Lemma
\ref{mia} below).

It is clear from the proofs given here for Neumann problem, that
one can deal with other situations such as Dirichlet, periodic or
mixed boundary conditions (see \cite{cavi} for scalar equations).
Systems like (\ref{s0}) have been considered also in \cite{clhi2}
and \cite{clhi}, where the matrix $A(x)$ is not necessarily
symmetric and with boundary conditions either of Dirichlet type or
of antiperiodic type. The authors establish sufficient conditions
for the positivity of the corresponding lower eigenvalue. These
conditions involve $L_1$ restrictions on the spectral radius of
some appropriate matrices which are calculated by using the matrix
$A(x).$ It is easy to check that, even in the scalar case, these
conditions are independent from classical $L_1-$Lyapunov
inequality (\ref{c1}) and therefore, for the ordinary case, they
are also independent from our results in this paper. Also, in a
series of papers, W. T. Reid (\cite{reid1}, \cite{reid2},
\cite{reid3}) made an extension of (\ref{c1}) for the Dirichlet
problem, but he always considered $p=1$ (see Remark
\ref{remark2107081} below).

In Section 3 we deal with elliptic systems of the form
\begin{equation}\label{spartialn} \Delta u(x) + A(x)u(x) = 0, \ x \in \Omega, \ \dis
\frac{\partial u(x)}{\partial n} = 0, \ x \in \partial \Omega
\end{equation} where $\Omega$ is a bounded and regular domain in
$\real^{N}$ and $\dis \frac{\partial}{\partial n}$ is the outer
normal derivative on $\partial \Omega.$ Here the relation between
$p$ and the dimension $N$ may be important (see Lemma \ref{miap}).
To our knowledge, there are no previous work on $L_p-$Lyapunov
inequalities for elliptic systems if $p \neq \infty$ (see
\cite{ba} and \cite{kawa}, section 5, for the case $p = \infty$).
Finally, we show some applications to nonlinear resonant problems.
In particular, and for Neumann boundary conditions, we obtain a
generalization for systems of equations of the main result given
in \cite{mawawi} where the author treated the scalar case and
where they use in the proof the duality method of Clarke and
Ekeland (see Theorem \ref{nolinealp} below).

\section{Ordinary boundary value problems}

This section will be concerned with boundary value problems of the
form (\ref{s1}). We begin with a preliminar lemma on
$L_\infty-$Lyapunov inequalities for (\ref{s1}), inspired from
\cite{la} and \cite{lasa}, where the authors studied periodic
boundary conditions. Our proof suggests the way to obtain optimal
$L_p-$Lyapunov inequalities for system (\ref{s1}) in the case $1
\leq p < \infty.$

\begin{lemma}\label{lpreliminar}
Let $A(\cdot)$ be a real $n\times n$ symmetric matrix valued
function with elements defined and continuous on $[0,L].$ Suppose
there exist diagonal matrix functions $P(x)$ and $Q(x)$ with
continuous respective entries $\delta_{kk} (x), \ 1 \leq k \leq n$
and $\mu_{kk} (x), \ 1 \leq k \leq n,$ and eigenvalues
$\lambda_{p(k)}, \ 1 \leq k \leq n,$ of the eigenvalue problem
(\ref{n2}) such that
\begin{equation}\label{eq1906081}
P(x) \leq A(x) \leq Q(x), \ \forall \  x \in [0,L]
\end{equation}
and
\begin{equation}\label{eq1906082}
\lambda_{p(k)} < \delta_{kk} (x) \leq \mu_{kk}(x) <
\lambda_{p(k)+1}, \ \forall \ x \in [0,L], \ 1 \leq k \leq n.
\end{equation}
Then there exists no nontrivial solution of (\ref{s1}).
\end{lemma}

\begin{proof}
Let us denote by $\sobo1$ the usual Sobolev space. If $u =
(u_1,\cdots,u_n) \in (H^{1}(0,L))^n,$ is a nontrivial solution of
(\ref{s1}), then
\begin{equation}\label{eq1906085} \displaystyle \int_{0}^{L}
<u'(x),v'(x)> \ dx = \displaystyle \int_{0}^{L}<A(x)u(x),v(x)>, \
dx \ \forall \ v \in \sobol,
\end{equation}
where $<\cdot,\cdot
>$ is the usual scalar product in $\real^n$.
The eigenvalues of (\ref{n2}) are given by $\lambda_j = \frac{j^2
\pi^2}{L^2},$ where $j$ is an arbitrary nonnegative integer
number. If $\varphi_j$ is the corresponding eigenfunction to
$\lambda_j,$ let us introduce the space $H = H_1 \times \cdots
\times H_k \times \cdots \times H_n,$ where for each $1 \leq k
\leq n,$ $H_k$ is the span of the eigenfunctions $\varphi_0,
\varphi_1, \cdots, \varphi_{p(k)}.$ It is trivial that we can
choose $\psi = (\psi_1,\cdots,\psi_n) \in H$ satisfying
\begin{equation}\label{eq1906084}u_k + \psi_k \in H_k^\bot, \ 1
\leq k \leq n. \end{equation} In fact
\begin{equation}\label{eq2006081}
\psi_k = \displaystyle \sum_{m=0}^{p(k)} c_m ^k \varphi_m, \ \ c_m
^k = \displaystyle - \frac{\intl u_k (x)\varphi_m (x) \ dx }{\intl
\varphi_m ^2  (x) \ dx}, \ 1 \leq k \leq n, \ 0 \leq m \leq p(k).
\end{equation}

The main ideas to get a contradiction with the fact that $u$ is a
nontrivial solution of (\ref{s1}) are the following two
inequalities. The first one is a consequence of the variational
characterization of the eigenvalues of (\ref{n2}). The second one
is a trivial consequence of the definition of the subspace $H_k.$
\begin{equation}\label{eq2006082}
\begin{array}{c}
\displaystyle \intl ((u_k + \psi_k)' (x))^2 \ dx \geq
\lambda_{p(k)+1} \displaystyle \intl ((u_k + \psi_k)(x))^2 \ dx, \
1 \leq k \leq n, \\ \\ \displaystyle \intl (\psi_k)' (x))^2 \ dx
\leq \lambda_{p(k)} \displaystyle \intl ((\psi_k)(x))^2 \ dx, \ 1
\leq k \leq n.
\end{array}
\end{equation}
Now, from (\ref{eq1906085}) we have
\begin{equation}\label{eq1906086}
\begin{array}{c}
\intl < (u+\psi)'(x),(u+\psi)'(x) > \ dx = \intl
<A(x)(u+\psi)(x),(u+\psi)(x)> \ dx + \\ \\
\intl < \psi'(x),\psi '(x)> \ dx - \intl < A(x)\psi (x),\psi(x)
> \ dx
\end{array}
\end{equation}
By using  (\ref{eq1906081}) and (\ref{eq1906082}) we deduce
$$
\begin{array}{c}
\intl < \psi'(x),\psi '(x)> \ dx - \intl < A(x)\psi (x),\psi(x)
> \ dx \leq \\ \\ \intl < \psi'(x),\psi '(x)> \ dx - \intl < P(x)\psi (x),\psi(x)
> \ dx \\ \\
= \displaystyle \sum_{k=1}^{n} \displaystyle \intl [(\psi_k '
(x))^2 - \delta_{kk}(x)(\psi_k (x))^2 ] \ dx  \leq \displaystyle
\sum_{k=1}^{n} \displaystyle \intl (\lambda_{p(k)}
-\delta_{kk}(x))(\psi_k (x))^2 \ dx \leq 0.
\end{array}
$$
Consequently,
\begin{equation}\label{eq1906086bis}
\intl < (u+\psi)'(x),(u+\psi)'(x) > \ dx \leq  \intl
<A(x)(u+\psi)(x),(u+\psi)(x)> \ dx.
\end{equation}
Also, from (\ref{eq1906081}), (\ref{eq1906084}), (\ref{eq2006082})
and (\ref{eq1906086bis}) we obtain
\begin{equation}\label{eq1906087}
\begin{array}{c}
\displaystyle \sum_{k=1}^{n} \intl \lambda_{p(k)+1} (u_k +
\psi_k)^2(x) \ dx \leq \displaystyle \sum_{k=1}^{n} \intl (u_k
+ \psi_k)'^{2} (x) \ dx = \\ \\
\intl < (u+\psi)'(x),(u+\psi)'(x)
> \ dx \leq  \intl
<A(x)(u+\psi)(x),(u+\psi)(x)> \ dx \leq \\ \\
\intl <Q(x)(u+\psi)(x),(u+\psi)(x)> \ dx = \displaystyle
\sum_{k=1}^{n} \intl \mu_{kk} (x) (u_k + \psi_k)^2(x) \ dx
\end{array}
\end{equation}
It follows, again from (\ref{eq1906082}), that
\begin{equation}\label{eq1906088}
u+\psi \equiv 0.
\end{equation}
But if $u + \psi \equiv 0,$ then $u = \phi =
(\phi_1,\cdots,\phi_n)$ for some nontrivial $\phi \in H.$
Therefore,
\begin{equation}\label{eq1906089}
\begin{array}{c}
\displaystyle \sum_{k=1}^{n} \intl \lambda_{p(k)} (\phi_k)^2(x) \
dx \geq \displaystyle \sum_{k=1}^{n} \intl (\phi_k)'^{2} (x) \ dx = \\ \\
\intl < \phi'(x),\phi'(x) \ dx =  \intl
<A(x)\phi(x),\phi(x)> \ dx \geq \\ \\
\intl <P(x)\phi(x),\phi(x)> \ dx = \displaystyle \sum_{k=1}^{n}
\intl \delta_{kk} (x) (\phi_k)^2(x) \ dx
\end{array}
\end{equation}
Now, (\ref{eq1906082}) implies that $u_k = \phi_k \equiv 0, \ 1
\leq k \leq n,$ which is a contradiction with the fact that $u$ is
nontrivial.
\end{proof}

\begin{remark}
It is clear from the previous proof that if the matrix functions
$P(x)$ and $Q(x)$ are constant functions $P$ and $Q$, then it is
not necessary to assume that they are, in addition, diagonal
matrices. In fact, to carry out the proof, it is sufficient to
assume that they are symmetric matrices and such that if
$\delta_k, \ 1 \leq k \leq n$ and $\mu_k, \ 1 \leq k \leq n$
denote the eigenvalues of $P$ and $Q$ respectively, then
\begin{equation}
\lambda_{p(k)} < \delta_{k} \leq \mu_{k} < \lambda_{p(k)+1}, \ 1
\leq k \leq n.
\end{equation}
\end{remark}

\skipaline

 We collect now some results which have been proved in
\cite{camovi}, section 2. Really, if we are treating with Lyapunov
inequalities for scalar ordinary problems and $1 \leq p < \infty,$
the constant $\beta_p$ defined in the next lemma, plays the same
role as $\beta_{\infty} = \lambda_1,$ in the $L_\infty-$Lyapunov
inequality (\ref{2806072})-(\ref{c2}).

\begin{lemma}\label{mia}(\cite{camovi})
If $1 \leq p \leq \infty$ is a given number, let us define the set
$X_p$ and the functional $I_p$ as
\begin{equation}\label{eq4}
\begin{array}{c}
X_{1} = \{ v \in \sobo1 : \max_{x \in [0,L]} v(x) + \min_{x \in
[0,L]} v(x) = 0 \}, \\ I_1: X_1 \setminus \{ 0 \} \rightarrow
\real, \ I_1(v) = \displaystyle \frac{\dis \intl v'^{2}}{\Vert v
\Vert_{\infty}^{2}}, \\ \\ X_{p} = \left \{ v \in \sobo1: \dis
\intl \vert v \vert
^{\frac{2}{p-1}} v = 0 \right \},\ \mbox{if} \ \ 1 < p < \infty, \\
I_p: X_p \setminus \{ 0 \} \rightarrow \real,I_{p}(v) = \dis
\frac{\dis \intl v'^{2}}{\left ( \dis \intl \vert v \vert
^{\frac{2p}{p-1}} \right )^{\frac{p-1}{p}}}, \ \mbox{if} \ \ 1 < p
< \infty, \\ \\ X_{\infty} = \{ v \in \sobo1: \ \dis \intl v = 0
\},
\\ I_\infty: X_\infty \setminus \{ 0 \} \rightarrow \real,\
I_{\infty}(v) = \displaystyle \frac{\dis \intl
v'^{2}}{\displaystyle \intl v^{2}}
\end{array}
\end{equation}
If
\begin{equation}\label{eq5}
\beta_{p} \equiv \min_{X_p \setminus \{ 0 \}} \ I_{p}, \ 1 \leq p
\leq \infty,
\end{equation}
and for some $p \in [1,\infty],$ function $a$ satisfies
(\ref{2806072}) and $\Vert a^+ \Vert_p < \beta_p,$ then (\ref{n1})
has only the trivial solution.
\end{lemma}

\begin{remark}\label{2409081}
It is possible to obtain an explicit expression for $\beta_p$, as
a function of $p$ and $L$ (see \cite{camovi}). In particular,
$\beta_1 = 4/L, \ \beta_\infty = \pi^2/L^2$ and $\beta_{1}$ is
attained in a function $v\in X_{1}\setminus \{ 0 \}$ if and only
there exists a nonzero constant $c$ such that $v(x) =
c(x-\frac{L}{2}), \forall x \in [0,L].$ Finally and in relation to
Lyapunov inequalities, the constant $\beta_p$ is optimal in  the
following sense (see \cite{camovi}): if
$$\Sigma_p = \{ a
\in L^{p}(0,L)\setminus \{0\}: \dis \intl a(x) \ dx \geq 0 \
\mbox{and} \ (\ref{n1}) \ \mbox{has nontrivial solutions} \ \} $$
then
$$ \beta_{1} \equiv \inf_{a \in \Sigma_1} \ \Vert a^+ \Vert_1, \
\beta_{p} \equiv \min_{a \in \Sigma_p} \ \Vert a^+ \Vert_p, \ 1 <
p \leq \infty.$$
\end{remark}

\skipaline

We return to system (\ref{s1}). From now on, we assume that the
matrix function $A(\cdot) \in \Lambda$ where $\Lambda$ is defined
as

\afirm{[$\Lambda$]}{10cm}{The set of real $n\times n$ symmetric
matrix valued function $A(\cdot)$, with continuous element
functions $a_{ij}(x), \ 1 \leq i,j \leq n, \ x \in [0,L],$ such
that (\ref{s1}) has not nontrivial constant solutions and
$$ \int_0^L <A(x)k,k>  dx \geq 0, \
\forall \ k \in \mathbb{R}^n. $$ }

The main result of this section is the following.

\begin{theorem}\label{t1}
Let $A(\cdot) \in \Lambda$ be such that there exist a diagonal
matrix $B(x)$ with continuous entries $b_{ii}(x),$ and $p_i \in
[1,\infty], \ 1 \leq i \leq n,$ satisfying
\begin{equation}\label{eq2306083}
\begin{array}{c}
A(x) \leq B(x), \ \forall \ x \in [0,L], \\ \\ \Vert b_{ii}^+
\Vert_{p_i} < \beta_{p_i}, \ \mbox{if} \ p_i \in (1,\infty], \ \
\Vert b_{ii}^+ \Vert_{p_i} \leq  \beta_{p_i}, \ \mbox{if} \ p_i
=1.
\end{array}
\end{equation}
Then, there exists no nontrivial solution of the vector boundary
value problem (\ref{s1}).
\end{theorem}
\begin{proof}
If $u \in \sobol$ is any nontrivial solution of (\ref{s1}), we
have
$$
\displaystyle \int_{0}^{L} <u'(x),v'(x)>  = \displaystyle
\int_{0}^{L}<A(x)u(x),v(x)>, \ \forall \ v \in \sobol.
$$
In particular, we have
\begin{equation}\label{eq3}
\begin{array}{c}
\displaystyle \intl <u'(x),u'(x)>  = \displaystyle \intl
<A(x)u(x),u(x)>, \ \\ \dis \intl <A(x)u(x),k> = \dis \intl
<A(x)k,u(x)> = 0, \ \forall \ k \in \real^n
\end{array}
\end{equation}
Therefore, for each $k \in \real^n, $ we have
$$
\begin{array}{c}
\dis \intl <(u(x)+k)',(u(x)+k)'> = \dis \intl <u'(x),u'(x)> = \\
\\
\displaystyle \intl <A(x)u(x),u(x)>  \leq \displaystyle \intl
<A(x)u(x),u(x)> + \dis \intl <A(x)u(x),k>+ \\ \\  \dis \intl
<A(x)k,u(x)>+ \dis \intl <A(x)k,k> =\\ \\ \dis \intl
<A(x)(u(x)+k),u(x)+k> \leq \dis \intl <B(x)(u(x)+k),u(x)+k>.
\end{array}
$$
If $u = (u_i), $ then for each $i, 1 \leq i \leq n,$ we choose
$k_i \in \real$ satisfying $u_i + k_i \in X_{p_i},$ the set
defined in Lemma \ref{mia}. By using previous inequality, Lemma
\ref{mia} and H{\"o}lder inequality, we obtain
\begin{equation}\label{eq8}
\begin{array}{c}
\dis \sum_{i=1}^n \beta_{p_i} \Vert (u_{i} + k_i)^2
\Vert_{\frac{p_i}{p_i -1}} \leq \dis \sum_{i=1}^n \intl (u_i
(x)+k_i)'^2 \leq \\ \dis \sum_{i=1}^n \intl b_{ii}^+ (x)(u_i
(x)+k_i)^2 \leq \dis \sum_{i=1}^n \Vert b_{ii}^+ \Vert_{p_i} \Vert
(u_{i} + k_i)^2 \Vert_{\frac{p_i}{p_i -1}},
\end{array}
\end{equation}
where
$$
\begin{array}{c}
\frac{p_i}{p_i -1} = \infty, \ \ \mbox{if} \ p_i = 1 \\ \\
\frac{p_i}{p_i -1} = 1, \ \ \mbox{if} \  p_i = \infty. \end{array}
$$

Therefore from (\ref{eq2306083}) we have
\begin{equation}\label{eq9}
\dis \sum_{i=1}^{n} (\beta_{p_i} - \Vert b_{ii}^+ \Vert_{p_i})
\Vert (u_{i} + k_i)^2 \Vert_{\frac{p_i}{p_i -1}} \leq 0.
\end{equation}
On the other hand, since $u$ is a nontrivial function, $u+k$ is
also a nontrivial function. Indeed, if $u+k$ is identically zero,
we deduce that (\ref{s1}) has the nontrivial and constant solution
$-k$ which is a contradiction with the hypothesis $A(\cdot) \in
\Lambda.$

Now, if $u+k$ is nontrivial, some component, say, $u_j + k_j$ is
nontrivial. If $p_j \in (1,\infty],$ then $(\beta_{p_j} - \Vert
b_{jj}^+ \Vert_{p_j}) \Vert (u_{j} + k_j)^2 \Vert_{\frac{p_j}{p_j
-1}}$ is strictly positive and from (\ref{eq2306083}), all the
other summands in (\ref{eq9}) are nonnegative. This is a
contradiction.

If $p_j =1,$ since $\beta_1$ is only attained in nontrivial
functions of the form $v(x) = c(x-\frac{L}{2}),$ and $v'(0) \neq
0,$ we have
$$\beta_{p_j} \Vert (u_{j} + k_j)^2 \Vert_{\frac{p_j}{p_j -1}} <
\intl (u_j (x)+k_j)'^2.$$

Then (\ref{eq8}) and (\ref{eq9}) are both strict inequalities and
this is again a contradiction.
\end{proof}

\begin{remark}\label{optimalidad}
Previous Theorem is optimal in the following sense. For any given
positive numbers $\gamma_i, \ 1 \leq i \leq n,$ such that at least
one of them, say $\gamma_j,$ satisfies
\begin{equation}\label{2optimalidad}
\gamma_j > \beta_{p_j}, \ \mbox{for some} \ p_j \in [1,\infty],
\end{equation} there exists a diagonal  $n\times n$ matrix
$A(\cdot) \in \Lambda$ with continuous entries $a_{ii}(x), \ 1
\leq i \leq n,$ satisfying $\Vert a_{ii}^+ \Vert_{p_i} <
\gamma_{i}, \ 1 \leq i \leq n$ and such that the boundary value
problem (\ref{s1}) has nontrivial solutions. To see this, if
$\gamma_j$ satisfies (\ref{2optimalidad}), then there exists some
continuous function $a(x), $ not identically zero, with $\intl a
(x) \ dx \geq 0, $ and $\Vert a^+ \Vert_{p_j} < \gamma_j,$ such
that the scalar problem
$$w''(x) + a(x)w(x) = 0, \ x \in (0,L), \ \ w'(0) = w'(L) = 0,
$$
has nontrivial solutions (see the remark after Lemma \ref{mia}).
Then, to get our purpose, it is sufficient to take $a_{jj}(x) =
a(x)$ and $a_{ii}(x) = \delta \in \real^+,$ if $i\neq j,$ with
$\delta$ sufficiently small.
\end{remark}

As an application of Theorem \ref{t1} we have the following
corollary.

\begin{corollary}\label{coro1}
Let $A(\cdot) \in \Lambda$ and, for each $x \in [0,L],$  let us
denote by $\rho (x)$ the spectral radius of the matrix $A(x).$ If
the function $\rho(\cdot)$ satisfies one of the following
conditions:
\begin{enumerate}
\item $ \Vert \rho^+ \Vert_{1} \leq \beta_{1},$%
\item There is some $p \in (1,\infty]$ such that $\Vert \rho^+
\Vert_{p} < \beta_{p}$,
\end{enumerate}
Then there exists no nontrivial solution of (\ref{s1}).
\end{corollary}
\begin{proof}
It is trivial, taking into account the previous Theorem and the
inequality
\begin{equation}\label{eqespectral}
A(x) \leq \rho (x) I_n, \ \forall \ x \in [0,L],
\end{equation}
where $I_n$ is the $n\times n$ identity matrix.
\end{proof}

\begin{remark}
The authors introduced in \cite{la} and \cite{lasa} similar
conditions for periodic problems and $p_i = \infty, \ 1 \leq i
\leq n$. Our method of proof, where we strongly use the
minimization problems considered in Lemma \ref{mia}, does possible
the consideration of the cases $p \in [1,\infty)$, which to the
best of our knowledge are new. In particular, if $p \in
[1,\infty),$ the function $\rho (x)$ may cross an arbitrary number
of eigenvalues of the problem (\ref{n2}). Also, by using our
methods one can deal with other boundary conditions and more
general second order equations (see, for the scalar case, Remark 5
in \cite{camovi} and Theorem 2.1 in \cite{cavi}).
\end{remark}

\begin{remark}\label{remark2107081}
In this remark we show some relations between previous Corollary
and some results contained in \cite{reid1}, \cite{reid2} and
\cite{reid3} for Dirichlet boundary conditions.

If $A(\cdot)$ satisfies \afirm{[H]}{10cm}{$A(x), \ x \in [0,L] $
is a continuous and positive semi-definite matrix function such
that $\det A(x) \neq 0 $ for some $x \in [0,L]$} (here $\det A(x)$
means the determinant of the matrix $A(x)$) and
\begin{equation}\label{eq2107081}
\intl \mbox{trace} \ A(x) \ dx \leq \beta_1,
\end{equation}
then there exists no nontrivial solution of (\ref{s1}). In fact,
taking into account that for each $x \in [0,L],$ $\rho (x)$ is an
eigenvalue of the matrix $A(x)$ and that in this case all the
eigenvalues of $A(x), \ \lambda_1(x), \cdots, \lambda_n (x),$ are
nonnegative, we have $\rho (x) \leq \sum_{i=1}^n \lambda_i (x) =
\mbox{trace} \ A(x)$ (see \cite{lang2} for this last relation).
Therefore, from (\ref{eq2107081}) we obtain
\begin{equation}\label{eq2107083} \Vert \rho^+ \Vert_{1} = \intl \rho(x) \
dx \leq \beta_{1}.
\end{equation}
\end{remark}
Previous remark shows that, if we want to have a criterion
implying that (\ref{s1}) has only the trivial solution, then
(\ref{eq2107083}) is better than (\ref{eq2107081}).

\skipaline

As in the scalar case, it may be seen that for Dirichlet boundary
conditions, hypothesis [H] is not necessary. However, for Neumann
boundary conditions, a restriction like [H] is natural (see Remark
4 and Remark 5 in \cite{camovi}).

\skipaline

In Corollary \ref{c2p} of the next section it is shown how, for
elliptic systems, we can obtain optimal conditions without the
help of the spectral radius of the matrix $A(x).$ Obviously that
Corollary is also applicable to ordinary problems as (\ref{s1}).

\section{Elliptic systems}

This section will be concerned with linear boundary value problems
of the form
\begin{equation}\label{s1p}
\Delta u(x)+ A(x)u(x) = 0,  \ x\in \Omega,  \ \ \frac{\partial
u(x)}{\partial n}=0, \  x\in \partial \Omega,
\end{equation}
Here $\Omega \subset \real^N, \ N \geq 2$ is a bounded and regular
domain, $\dis \frac{\partial}{\partial n}$ is the outer normal
derivative on $\partial \Omega$ and $A \in \Lambda_{*}$, where
$\Lambda_{*}$ is defined as
\afirm{[$\Lambda_{*}$]}{10cm}{The set
of real $n\times n$ symmetric matrix valued function $A(\cdot)$,
with continuous element functions $a_{ij}(x), \ 1 \leq i,j \leq n,
\ x \in \comega,$ such that (\ref{s1p}) has not nontrivial
constant solutions and
\begin{equation}\label{04061} \intomega <A(x)k,k>  dx \geq 0, \ \forall \ k \in
\mathbb{R}^n. \end{equation} }  In (\ref{s1p}), $u\in
(H^{1}(\Omega))^n,$ the usual Sobolev space.

As in the ordinary case, we now collect some results which have
been proved in \cite{camovi2}.

\begin{lemma}\label{miap}(\cite{camovi2})
If $1 \leq \frac{N}{2} < p \leq \infty$ is a given number, let us
define the set $X_p$ and the functional $I_p$ as
\begin{equation}\label{eq4p}
\begin{array}{c}
X_{p} = \left \{ v \in H^1(\Omega): \dis \int_{\Omega} \vert v
\vert ^{\frac{2}{p-1}} v = 0 \right \}, \ \ \mbox{if} \ \
\frac{N}{2}
< p < \infty, \\ \\

I_p: X_p \setminus \{ 0 \} \rightarrow \real, \ I_{p}(v) = \dis
\frac{\dis \int_{\Omega} \vert \nabla v\vert^2}{\left ( \dis
\int_{\Omega} \vert v \vert ^{\frac{2p}{p-1}} \right
)^{\frac{p-1}{p}}}, \ \mbox{if} \ \ \frac{N}{2} < p < \infty, \\ \\
X_{\infty} = \{ v \in H^1(\Omega): \ \dis \intomega v = 0 \}, \\
I_\infty: X_\infty \setminus \{ 0 \} \rightarrow \real,\
I_{\infty}(v) = \displaystyle \frac{\dis \intomega \vert \nabla v
\vert ^2}{\displaystyle \intomega v^{2}}
\end{array}
\end{equation}
If
\begin{equation}\label{eq5p}
\beta_{p} \equiv \min_{X_p \setminus \{ 0 \}} \ I_{p},\
\frac{N}{2} < p \leq \infty,
\end{equation}
and a given function $a$ satisfies
\begin{equation}\label{eq2306081}
a \in L^p (\Omega,\real) \setminus \{ 0 \}, \ \int_{\Omega} a \geq
0, \ \Vert a^+ \Vert_{p} < \beta_p,
\end{equation}
then the scalar problem
\begin{equation}\label{eq2306082}
\Delta u(x)+ a(x)u(x) = 0,  \ x\in \Omega,  \ \ \frac{\partial
u(x)}{\partial n}=0, \  x\in \partial \Omega,
\end{equation}
 has only the trivial solution.

\end{lemma}

\begin{remark}\label{rp1}
As in the ordinary case, $\beta_\infty = \lambda_1,$ the first
strictly positive eigenvalue of the Neumann eigenvalue problem in
the domain $\Omega$. Consequently, it seems difficult to obtain
explicit expressions for $\beta_{p},$ as a function of $p,\Omega$
and $N,$ at least for general domains. Finally, the constant
$\beta_p$ is optimal in  the following sense: if $\frac{N}{2} < p
\leq \infty$ and
$$\Sigma_p ^* = \{ a
\in L^{p}(\Omega)\setminus \{0\}: \dis \intomega a(x) \ dx \geq 0
\ \mbox{and} \ (\ref{s1p}) \ \mbox{has nontrivial solutions} \ \}
$$ then
$$ \beta_{p} \equiv \min_{a \in \Sigma_p ^*} \ \Vert a^+ \Vert_p, \ N/2 <
p \leq \infty.$$
\end{remark}

Next result may be proved by using the same ideas as in Theorem
\ref{t1}.

\begin{theorem}\label{t1p}
Let $A(\cdot) \in \Lambda_{*}$ be such that there exist a diagonal
matrix $B(x)$ with continuous entries $b_{ii}(x),$ and numbers
$p_i \in (N/2,\infty], \ 1 \leq i \leq n,$ which fulfil
\begin{equation}\label{eq7p}
\begin{array}{c}
A(x) \leq B(x), \ \forall \ x  \in \ \comega
\\ \\ \Vert b_{ii}^+
\Vert_{p_i} < \beta_{p_i}, \ 1 \leq i \leq n.
\end{array}
\end{equation}
Then, there exists no nontrivial solution of the vector boundary
value problem (\ref{s1p}).
\end{theorem}

\begin{remark}\label{optimalidadp}
As in the ordinary case, the previous Theorem is optimal in the
sense of Remark \ref{optimalidad} (see Theorem 2.1 in
\cite{camovi2}). Moreover, by using the previous Theorem, it is
possible to obtain a corollary similar to corollary \ref{coro1},
which involves the spectral radius $\rho (x)$ of the matrix $A(x)$
and the norm $\Vert \rho^+ \Vert_p.$ The unique difference with
the ordinary case is that, for elliptic systems, $p \in
(N/2,\infty].$
\end{remark}
 In the next Corollary and in order to show how our Theorem
\ref{t1p} can be used without the help of the spectral radius of
the matrix $A(x),$ we consider the case of a system with two
equations.

\begin{corollary}\label{c2p} Let the matrix $A(x)$ be given by
\begin{equation}\label{eq10}
A(x) = \left (
\begin{array}{cc} a_{11}(x) & a_{12}(x) \\
a_{12}(x) & a_{22}(x) \end{array} \right )
\end{equation}
where \afirm{[{\bf H1}]}{10cm} {
$$
\begin{array}{c} a_{ij} \in C(\comega), \ 1 \leq i,j \leq 2, \\ \\ a_{11}(x) \geq 0, \
a_{22}(x) \geq 0, \ a_{11}(x)a_{22}(x) \geq a_{12}^2 (x), \
\forall \ x \in \comega, \\ \\ \det \ A(x) \neq 0, \ \mbox{for
some}\ x \in \comega. \end{array}
$$ }
In addition, let us assume that there exist $p_1, p_2 \in
(N/2,\infty]$ such that
\begin{equation}\label{eq11}
\Vert a_{11}\Vert_{p_1} < \beta_{p_1}, \ \ \Vert a_{22} +
\displaystyle \frac{a_{12}^2}{\beta_{p_1} - \Vert a_{11}
\Vert_{p_1}} \Vert _{p_2} < \beta_{p_2}.
\end{equation}
Then the unique solution of (\ref{s1p}) is the trivial one.
\end{corollary}

\begin{proof}
It is trivial to see that [{\bf H1}] implies that the eigenvalues
of the matrix $A(x)$ are both nonnegative, which implies that
$A(x)$ is positive semi-definite. Also, since  $\det \ A(x) \neq
0, \ \mbox{for some}\ x \in \comega,$ (\ref{s1p}) has not
nontrivial constant solutions. Therefore, $A(\cdot) \in
\Lambda_{*}.$ Moreover, it is easy to check that for a given
diagonal matrix $B(x),$ with continuous entries $b_{ii}(x), \ 1
\leq i \leq 2,$ the relation
\begin{equation}\label{eq2306084}
A(x) \leq B(x), \ \forall \ x \in \comega
\end{equation}
is satisfied if and only if $\forall \ x \in \comega,$ we have
\begin{equation}\label{eq12}
\begin{array}{c}
\ b_{11}(x) \geq a_{11}(x), \ b_{22}(x) \geq  a_{22}(x),  \\
(b_{11}(x)- a_{11}(x))(b_{22}(x)- a_{22}(x)) \geq a_{12}^2 (x).
\end{array}
\end{equation}
In our case, if we choose
\begin{equation}\label{eq13}
b_{11}(x) = a_{11}(x) + \gamma, \ b_{22}(x) = a_{22}(x) +
\displaystyle \frac{a_{12}^2 (x)}{\gamma}
\end{equation}
where $\gamma$ is any constant such that
\begin{equation}
\begin{array}{c}
0 < \gamma < \beta_{p_1} - \Vert a_{11} \Vert_{p_1}, \\
\left ( \frac{1}{\gamma} - \frac{1}{\beta_{p_1} - \Vert a_{11}
\Vert _{p_1}} \right ) \Vert a_{12}^2 \Vert_{p_2} < \beta_{p_2} -
\Vert a_{22} + \displaystyle \frac{a_{12}^2}{\beta_{p_1} - \Vert
a_{11} \Vert_{p_1}} \Vert _{p_2}
\end{array}
\end{equation}
then all conditions of Theorem \ref{t1p} are fulfilled and
consequently (\ref{s1p}) has only the trivial solution.
\end{proof}

\begin{remark}\label{r2}
Previous corollary may be seen as a perturbation result in the
following sense: let us assume that we have an uncoupled system of
the type
\begin{equation}\label{s2p}
\begin{array}{c}
\Delta u_1  (x) + a_{11}(x)u_1 (x) = 0, \ x \in \Omega; \ \
\frac{\partial u_1 (x)}{\partial n}=0\,  \ x\in \partial \Omega,

\\ \\  \Delta u_2  (x) + a_{22}(x)u_2 (x) = 0, \ x \in \Omega; \ \
\frac{\partial u_2 (x)}{\partial n}=0\,  \ x\in \partial \Omega,
\end{array}
\end{equation}
where \begin{equation}\label{eq14}
\begin{array}{c}
a_{ii} \in C(\comega), \ 1 \leq i \leq 2, \  a_{11}(x) \geq \delta
> 0, \ a_{22}(x) \geq \delta , \ \forall \ x \in \comega.
\\ \\
\exists \ p_1, p_2 \in (N/2,\infty]:\  \Vert a_{11}\Vert_{p_1} <
\beta_{p_1}, \ \ \Vert a_{22}  \Vert _{p_2} < \beta_{p_2}.
\end{array}
\end{equation}
Then it is clear from the scalar results (see Remark \ref{rp1})
that the unique solution of (\ref{s2p}) is the trivial one (see
Corollary 6.1 in \cite{camovi2}). Now, we can use Corollary
\ref{c2p} to ensure the permanence of the uniqueness property
(with respect to the existence of solutions) of the coupled system
(\ref{s1p}), for any function $a_{12} \in C(\comega)$ with
$L^\infty-$norm sufficiently small. Here we have considered that
the functions $a_{ii}(x), \ 1 \leq i \leq 2,$ are fixed and that
the uncoupled system is perturbed by the function $a_{12}(x).$ But
it is clear that we may consider, for example, $a_{11}(x),
a_{12}(x)$ fixed and $a_{22}(x)$ as the perturbation. Some of
these results may be generalized to systems with $n$ equations.
For example, if we have an uncoupled system of the type
\begin{equation}\label{biss2p}
\Delta u_i  (x) + a_{ii}(x)u_i (x) = 0, \ x \in \Omega; \ \
\frac{\partial u_i (x)}{\partial n}=0\,  \ x\in \partial \Omega, \
1 \leq i \leq n
\end{equation}
where \begin{equation}\label{eq14bis}
\begin{array}{c}
a_{ii} \in C(\comega), \ 1 \leq i \leq n, \  a_{ii}(x) \geq \delta
> 0, \ 1 \leq i \leq n, \ \forall \ x \in \comega,
\\ \\
\exists \ p_i, \in (N/2,\infty]:\  \Vert a_{ii}\Vert_{p_i} <
\beta_{p_i}, \ 1 \leq i \leq n,
\end{array}
\end{equation}
then we can use Theorem \ref{t1p} to ensure the permanence of the
uniqueness property (with respect to the existence of solutions)
of the coupled system (\ref{s1p}), for any functions $a_{ij} =
a_{ji} \in C(\comega), \ 1 \leq i \neq j \leq n$ with
$L^\infty-$norm sufficiently small. The proof is similar to the
case of two equations and it is based on Theorem \ref{t1p}. The
unique difference is that now, the matrix $B(x)$ is given by
$b_{ii}(x) = a_{ii}(x) + \varepsilon, \ 1 \leq i \leq n $ with
$\varepsilon$ sufficiently small. It is easily deduced that if the
$L^\infty-$norm of the functions $a_{ij}= a_{ji}, \ 1 \leq i \neq
j \leq n$ are sufficiently small,  then the matrix $B(x) - A(x)$
is positive definite for all $x \in \comega.$
\end{remark}

Next we give some new results on the existence and uniqueness of
solutions of nonlinear resonant problems. We prefer to deal with
systems of P.D.E. (similar results can be proved for ordinary
differential systems; in this last case it is possible to choose
the constants $p_i \in [1,\infty], \ 1 \leq i \leq n $). In
particular, next Theorem is a generalization, for systems of
equations, of the main result given in \cite{mawawi} for the
Neumann problem. Moreover, it is a generalization (at the two
first eigenvalues of (\ref{n2})) of some results given in
\cite{ba} and \cite{kawa} where the authors take all the constants
$p_i = \infty, \ 1 \leq i \leq n$.

In the proof, the basic idea is to combine the results obtained in
the linear case with Schauder's fixed point theorem.

\bigskip

\begin{theorem}\label{nolinealp}
Let $\Omega \subset \real^N$ ($N\geq 2$) be a bounded and regular
domain and $G : \comega \times \real^n \rightarrow \real, \ (x,u)
\rightarrow G(x,u) $ satisfying:
\begin{enumerate}
\item
\begin{enumerate}
\item $u \rightarrow G(x,u)$ is of class $C^2(\real^n,\real)$ for
every $x \in \comega.$ \item $x \rightarrow G(x,u)$ is continuous
on $\comega$ for every $u \in \real^n.$
\end{enumerate}
\item %
There exist continuous matrix functions $A(\cdot),$ $B(\cdot),$
with $B(x)$ diagonal and with entries $b_{ii}(x),$ and $p_i \in
(N/2,\infty] \ 1 \leq i \leq n,$ such that
\begin{equation}\label{eq15}
\left.
\begin{array}{c}
A(x) \leq G_{uu}(x,u) \leq B(x) \ \mbox{in} \ \comega \times
\real^n, \\ \\ \Vert b_{ii}^+ \Vert_{p_i} < \beta_{p_i}, \ 1 \leq
i \leq n, \\ \\
\intomega <A(x)k,k>  dx > 0, \ \forall \ k \in \mathbb{R}^n
\setminus \{ 0 \}.
\end{array}
\right \}
\end{equation}
\end{enumerate}
Then system
\begin{equation}\label{eq17}
\left.
\begin{array}{c}
\Delta u(x)+ G_u (x,u(x)) = 0,  \ x\in \Omega,  \\
\frac{\partial u(x)}{\partial n}=0, \  \ x\in \partial \Omega,
\end{array}
\right \}
\end{equation}
has a unique solution.
\end{theorem}

\begin{proof}  We first prove uniqueness. Let $v$ and $w$ be two
solutions of (\ref{eq17}). Then, the function $u = v-w$ is a
solution of the problem
\begin{equation}\label{eq18}
\Delta u(x) + C(x)u(x)= 0 , \ x \in \Omega, \ \frac{\partial
u}{\partial n}= 0, \ x \in \partial \Omega
\end{equation}
where  $C(x) = \dis \int_{0}^{1} G_{uu}(x,w(x)+\theta u(x)) \ d
\theta$ (see \cite{lang}, p. 103, for the mean value theorem for
the vectorial function $G_u (x,u)$). Hence $A(x)\leq C(x)\leq
B(x)$ and we deduce that $C(x)$ satisfies all the hypotheses of
Theorem \ref{t1p}. Consequently, $u \equiv 0.$
\newline
Next we prove existence. First, we write (\ref{eq17}) in the
equivalent form
\begin{equation}\label{eq19} \left.
\begin{array}{cc}
\Delta u(x) + D(x,u(x))u(x) + G_u(x,0) = 0, & \text{ in } \Omega, \\
 \frac{\partial u}{\partial n}= 0, & \text{ on } \partial \Omega
\end{array}\right\}
\end{equation}
where the  function $D: \comega \times \real^n \rightarrow
\mathcal{M}(\real)$ is defined by $D(x,z) = \dis \int_{0}^{1}
G_{uu}(x,\theta z) \ d \theta.$ Here $\mathcal{M}(\real)$ denotes
the set of real $n \times n$ matrices. Let $X=(C(\comega))^n$ be
with the uniform norm, i.e., if $y(\cdot) = (y^1
(\cdot),\cdots,y^n (\cdot)) \in X, $ then $\Vert y \Vert_X =
\displaystyle \sum_{k=1}^n \Vert y^k (\cdot) \Vert_{\infty}.$
Since
\begin{equation}\label{eq2306085}
A(x) \leq D(x,z) \leq B(x), \ \forall \ (x,z) \in \comega \times
\real^n,
\end{equation}
we can apply Theorem \ref{t1p} in order to have a well defined
operator $T: X \rightarrow X,$ by $Ty = u_{y}$, being $u_{y}$ the
unique solution of the linear problem
\begin{equation}\label{eq20}\left.
\begin{array}{cc}
\Delta u(x)+ D(x,y(x))u(x) + G_u (x,0)= 0, & \text{ in } \Omega, \\
 \frac{\partial u}{\partial n}= 0, & \text{ on } \partial \Omega.
\end{array}\right\}
\end{equation}
We will show that $T$ is completely continuous and that $T(X)$ is
bounded. The Schauder's fixed point theorem provides a fixed point
for $T$ which is a solution of (\ref{eq17}).

The fact that $T$ is completely continuous is a consequence of the
compact embedding of the Sobolev space $W^{2,q}(\Omega)\subset
C(\comega)$ for $q$ sufficiently large. It remains to prove that
$T(X)$ is bounded. Suppose, contrary to our claim, that $T(X)$ is
not bounded. In this case, there would exist a sequence $\{ y_{n}
\} \subset X$ such that $\Vert u_{y_{n}}\Vert _{X} \rightarrow
\infty.$ From (\ref{eq2306085}), and passing to a subsequence if
necessary, we may assume that, for each $1 \leq i,j \leq n,$ the
sequence of functions $\{ D_{ij}(\cdot,y_{n}(\cdot))\}$ is weakly
convergent in $L^{p}(\Omega)$ to a function $E_{ij}(\cdot)$ and
such that if $E(x) = (E_{ij}(x)),$ then $A(x) \leq E(x) \leq
B(x)$, a.e. in $\Omega, $ (\cite{lema}, page 157).

If $z_{n} \equiv \dis \frac{u_{y_{n}}}{\Vert u_{y_{n}} \Vert
_{X}},$ passing to a subsequence if necessary, we may assume that
$z_{n} \rightarrow z_{0}$ strongly in $X$ (we have used again the
compact embedding $W^{2,q}(\Omega)\subset C(\comega)$), where
$z_{0}$ is a nonzero vectorial function satisfying
\begin{equation}\label{eq34}\left .
\begin{array}{cc}
 \Delta z_{0}(x)+ E(x)z_{0}(x)= 0,& \mbox{ in } \Omega, \\
\frac{\partial z_{0}}{\partial n}= 0, & \mbox{ on } \partial
\Omega
\end{array}\right \}
\end{equation}
This is a contradiction with Theorem \ref{t1p}.
\end{proof}

\end{document}